\documentclass{article}
\usepackage{euscript,amsfonts,amssymb,amsmath,multicol,amsthm,slashbox,longtable,amscd,url}
\usepackage[all]{xy}
\usepackage{graphicx,epsfig}
\author{Alexei G. Gorinov}
\title{Division theorems for the rational cohomology of discriminant complements and applications to automorphism groups in finite characteristic}
\date{}

\newcommand{\PP}{{\mathbb P}}
\newcommand{\Z}{{\mathbb Z}}
\newcommand{\C}{{\mathbb C}}

\newcommand{\Q}{{\mathbb Q}}
\newcommand{\dd}{\underline{d}}
\newcommand{\aaa}{\mathbf{a}}

\newcommand{\cc}{\mathbf{c}}

\newcommand{\m}{\mathbf{m}}

\newcommand{\LCM}{\mathrm{LCM}}

\newtheorem{theorem}{Theorem}

\newtheorem{lemma}{Lemma}\textwidth=14truecm
\newtheorem{Prop}{Proposition}

\newtheorem{conjecture}{Conjecture}
\DeclareMathOperator\Sing{Sing}

\newcommand{\SL}{\mathrm{SL}}
\newcommand{\PGL}{\mathrm{PGL}}
\newcommand{\GL}{\mathrm{GL}}

\begin{document}
\maketitle
\begin{abstract}
In this note we extend some of the results of a previous paper \url{arXiv:math/0511593} to algebraically closed fields of finite characteristic. In particular, we show that there is an explicit expression in $n$ and $d$ which is divisible by the prime to $p$ part of the order of the the automorphism group of a smooth degree $d$ hypersurface $\subset \mathbb{P}^n_k$ for $k$ an algebraically closed field  of characteristic $p$.
\end{abstract}
\section{Introduction and main results}\label{intro}

To state our main results let us first recall some of the notation of \cite{division}.
\label{nkdd} We fix integers $n$ and $k$ be  satisfying $1\leq k\leq n+1$ and let $F$ be an algebraically closed field. In the sequel $n$ will be the dimension of the ambient projective space and $k$ the codimension of a smooth complete intersection.

Set $\dd=(d_1,\ldots,d_k)$ to be a collection
of integers $\geq 2$. Let \label{PiSigma}$\Pi_{\dd,n}(F)$ be the $F$-vector
space of all $k$-tuples $(f_1,\ldots ,f_k)$, where $f_i,i=1,\ldots, k$, is a homogeneous polynomial
in $n+1$ variables of degree $d_i$ with coefficients in $F$. For every $(f_1,\ldots ,f_k)\in\Pi_{\dd,n}(F)$
let \label{SING}$\Sing (f_1,\ldots ,f_k)$ be the projectivisation of the set of all $x\in F^{n+1}\setminus\{ 0\}$ such that
\begin{itemize}
\item $f_i(x)=0,i=1,\ldots,k$,
\item the gradients of $f_i,i=1,\ldots,k$ at $x$ are linearly dependent.
\end{itemize}
Set $\Sigma_{\dd,n}(F)$ to be the
subset of $\Pi_{\dd,n}$ consisting of all $(f_1,\ldots ,f_k)$ such that
$\Sing (f_1,\ldots ,f_k)\neq\varnothing$.
%For any $(f_1,\ldots,f_k)\in\Pi_{\dd,n}$ set $\langle f_1\ldots,f_k\rangle$ to be the ideal generated by $f_1,\ldots ,f_k$.

Recall that a complete intersection $\subset \mathbb{P}^n (F)$ given by $f_1=\cdots =f_k=0$ is singular if and only if $(f_1,\ldots, f_k)\in\Sigma_{\dd,n}(F)$ (see e.g. Hartshorne \cite[exercise 5.8, chapter 1]{hartshorne}).

In the sequel, when $F=\C$ we will simply write $\Pi_{\dd,n}$ and $\Sigma_{\dd,n}$ instead of $\Pi_{\dd,n}(F)$, respectively $\Sigma_{\dd,n}(F)$.

The group $\GL_{n+1}(\C)$ acts on $\Pi_{\dd,n}$ in a natural way:
\begin{equation}\label{firstaction}\Pi_{\dd,n}\times\GL_{n+1}(\C)\ni ((f_1,\ldots,f_k),A)\mapsto (f_1\circ A,\ldots, f_k\circ A);\end{equation}
this action preserves $\Sigma_{\dd,n}$ (and hence, $\Pi_{\dd,n}\setminus\Sigma_{\dd,n}$).
%We will write $\Pi_{\dd,n}$ and
%$\Sigma_{\dd,n}$ instead of $\Pi_{\dd,n}(\C)$ and $\Sigma_{\dd,n}(\C)$ respectively.
%In the sequel we work almost entirely over $\C$, however some of the results can be extended to the case of abritrary $\kk$ completely
%for free, which is why we introduced the notation in a more general setting.

The following theorem was proved in \cite{division}:

\begin{theorem}\label{main}
Suppose $\dd\neq (2)$. Then the geometric quotient of $\Pi_{\dd,n}\setminus\Sigma_{\dd,n}$ by $\GL_{n+1}(\C)$ exists, and
the Leray spectral sequence with $\Q$-coefficients of the corresponding quotient map degenerates
at the second term.
\end{theorem}
This theorem generalises a result of J. Steenbrink and C. Peters in the case $k=1$ \cite{stepet}. 
%Our general strategy will be the same as in \cite{stepet}; however, the details will be different and more elementary (or so we hope).

The proof is based on the Leray-Hirsch principle: it suffices to construct global cohomology classes on $\Pi_{\dd,n}\setminus\Sigma_{\dd,n}$ such that their pullbacks under any orbit map generate the cohomology of the group $\GL_{n+1}(\C)$ (as a topological space).
We construct such classes by taking the Alexander duals of certain
subvarieties of~$\Sigma_{\dd,n}$. (These subvarieties are formed by all elements of $\Pi_{\dd,n}$ such that the corresponding complete intersections have singularities at each point of some projective subspace and maybe elsewhere, see section 2 of \cite{division}.) The resulting classes (they were denoted  $\aaa_i^{\dd,n}$ in \cite{division}) are manifestly Tate.

%\begin{theorem}\label{maingeneral}
%Let $W_{\dd,n}$ be the space of all smooth complete intersection of multidegree $\dd$ in $\mathbb{P}^n(\C)$ (we consider $W_{\dd,n}$ as a subvariety of the Chow scheme of degree $d_1\cdot\cdots\cdot d_k$ subvarieties of $\mathbb{P}^n(\C)$). Then the stack $[W_{\dd,n}/\PGL_{n+1}(\C)]$ is separated and Deligne-Mumford. Let $M_{\dd,n}$ be its coarse moduli space. We have an isomorphism $$H^*(W_{\dd,n},\Q)\cong H^*(\PGL_{n+1}(\C),\Q)\otimes H^*(M_{\dd,n},\Q).$$
%\end{theorem}

%\medskip

The following statements are by-products of the proof.

\newcounter{vb}
\begin{theorem}\label{upperboundvector}
Let $d$ be an integer $>2$.% and let $L$ be an algebraically closed field of characteristic 0. Then the
The order of the subgroup of $\GL_{n+1}(\C)$ consisting of the transformations that fix $f\in\Pi_{(d),n}\setminus\Sigma_{(d),n}$ divides
$$\prod_{i=0}^n((-1)^{n-i}+(d-1)^{n-i+1})(d-1)^i.$$
\end{theorem}
\setcounter{vb}{\value{theorem}}

\newcounter{pb}
\begin{theorem}\label{upboundprojective}
The order of the subgroup $\PGL_{n+1}(\C),n\geq 1,$ consisting of the transformations that preserve a smooth hypersurface of degree $d>2$ divides
\begin{equation}\label{bound}
\frac{1}{n+1}\prod_{i=0}^{n-1}\frac{1}{C^i_{n+1}}((-1)^{n-i}+(d-1)^{n-i+1})\LCM (C^i_{n+1}(d-1)^i,(n+1)(d-1)^n).
\end{equation}
%In particular, if $d\neq n+1,n\geq 3$, then the order of the automorphism group of $X$ divides~(\ref{bound}).
\end{theorem}
\setcounter{pb}{\value{theorem}}

Here $\LCM$ stands for the least common multiple and $C_{n+1}^i$ are the binomial coefficients: $C_{n+1}^i=\frac{i!(n+1-i)!}{(n+1)!}$. 

The purpose of the present paper is to extend these results to finite characteristic.

If $n\geq 3,d\geq 3$ and $(d,n)\neq(4,3)$, then any automorphism of a smooth hypersurface of degree $d$ in $\PP^n(\C)$ is known \cite[theorem 2]{mm} to be the restriction of a projective
transformation, so in these
cases theorem \ref{upboundprojective} implies that the order of the full automorphism group divides~(\ref{bound}). We later prove this for arbitrary complete intersections in arbitrary charactiristic, see lemma \ref{projautisaut}.

The proofs of the above results in the complex case are based on the following observation.

\newtheorem*{observation}{Observation}

\begin{observation}\label{obser}
Suppose a connected Lie group $G$ acts on a pathwise connected topological space $X$ and let $m$ be the dimension of a maximal compact subgroup $G^c\subset G$. Then (since $G$ as a topological space can be contracted onto $G^c$) we have $H^m(G,\Z)\cong \Z$. Suppose that for some choice of $x\in X$ the integral cohomology map induced by the orbit map $G\ni g\mapsto g\cdot x\in X$ contains an element an element that spans the subgroup of index $a\in\Z$. Then the order of any finite subgroup of any stabiliser divides $a$.
\end{observation}

%This will be proved in section \ref{proofobs}.

The variety $\Sigma_{\dd,n}$ is in fact defined over $\Z$ (see lemma \ref{definedoverz}) so it seems natural to ask whether there are analogues of theorems \ref{main}, \ref{upperboundvector} and \ref{upboundprojective} in finite characteristic. It turns out that the answer to the latter two questions is positive but the analogues of theorems \ref{upperboundvector} and \ref{upboundprojective} give information on the prime to the characteristic part of the automorphism group only. The answer to the first one is likely to be positive as well but at the moment we are able to prove it only for stack cohomology, see lemma \ref{mainisoforstacks}.

Let $K$ be an algebraically closed field of characteristic $p$ and let $\ell\neq p$ be a prime. We write $\Pi_{\dd,n,K}$ to denote the space of $k$-tuples $(f_1,\ldots,f_k)$ where each $f_i$ is a homogeneous polynomial of degree $d_i$ with coefficients in $K$.

Take an irreducible polynomial $\triangle_{\dd,n}$ with integer coefficients that defines $\Sigma_{\dd,n}$; assume the greatest common divisor of the coefficients is 1.%\footnote{In fact, one can write down explicit defining equations for the variety $V$ from the proof of lemma \ref{definedoverz}. One can then use elimination theory to find $\triangle_{\dd,n}$ explicitly, cf. the book \cite{gkz} by I.~M.~Gelfand, M.~M.~Kapranov and A.~V.~Zelevinsky. This however will not be attempted in this paper.}.
Let $\Sigma_{\dd,n,K}$ be the subvariety of $\Pi_{\dd,n,K}$ given by $\triangle_{\dd,n}=0$. For a variety $X$ over $K$ and a finite abelian group $A$ of order prime to $p$ we denote the \'etale cohomology of $X$ with coefficients in $A$ by $H^*(X,A)$; as usual we write $H^*(X,\Z_\ell)$ for $\varprojlim H^*(X,\Z/\ell^i\Z)$ and $H^*(X,\Q_\ell)$ for $H^*(X,\Z_\ell)\otimes \Q_\ell$. I.e. we do not use the {\it \'et} subscript for the \'etale cohomology.% in finite characteristic. %(but we'll eventually use it when comparing the \'etale and singular cohomology in characteristic 0). 

%\begin{theorem}\label{lerayhirschcharp}
%We have a multiplicative isomorphism $$H^*(\Pi_{\dd,n,K}\setminus\Sigma_{\dd,n,K},\Q_\ell)\cong H^*(\GL_{n+1}(K),\Q_\ell)\otimes H^*((\Pi_{\dd,n,K}\setminus\Sigma_{\dd,n,K})/\GL_{n+1}(K),\Q_\ell)$$

%\end{theorem}

\begin{theorem}\label{upperboundvectorcharp}
Let $d$ be an integer $>2$. Then the prime to $p$ part of the
order of the subgroup of $\GL_{n+1}(K)$ consisting of the transformations that fix $f\in\Pi_{(d),n,K}\setminus\Sigma_{(d),n,K}$ divides
$$\prod_{i=0}^n((-1)^{n-i}+(d-1)^{n-i+1})(d-1)^i.$$
\end{theorem}

\begin{theorem}\label{upboundprojectivecharp}
1. The prime to $p$ part of the order of the subgroup $\PGL_{n+1}(K),n\geq 1,$ consisting of the transformations that preserve a smooth degree $d>2$ hypersurface of $\PP_n(K)$ divides (\ref{bound}).

2. If $n>2$ and $(n,d)\neq (3,4)$, then the projective automorphism group of a smooth degree $d$ projective hypersurface coincides with the full automorphinm group. \end{theorem}

Let us write explicitly the expression (\ref{bound}) for curves in $\PP^2$, surfaces in $\PP^3$ and threefolds in $\PP^4$
\begin{equation}\label{krivye}d^2(d-1)^4(d^2-3d+3)(d-2),\end{equation}
$$\frac{1}{3}d^3(d-1)^8(d^3-4d^2+6d-4)(d^2-3d+3)(d-2)\LCM(3,2(d-1))$$
and
\begin{multline*}\frac{1}{4}d^4(d-1)^{13}(d^4-5d^3+10d^2-10d+5)(d^3-4d^2+6d-4)(d^2-3d+3)(d-2)\\
\cdot\LCM(2,(d-1)^2)\LCM(2,d-1)\end{multline*}

The author is grateful to Bas Edixhoven for useful conversations and correspondence. This is a preliminary version of the paper. In particular, the following statement seems likely to be true but at this stage remains a conjecture.

\begin{conjecture}\label{lerayhirschcharp}
We have a multiplicative isomorphism $$H^*(\Pi_{\dd,n,K}\setminus\Sigma_{\dd,n,K},\Q_\ell)\cong H^*(\GL_{n+1}(K),\Q_\ell)\otimes H^*((\Pi_{\dd,n,K}\setminus\Sigma_{\dd,n,K})/\GL_{n+1}(K),\Q_\ell)$$

\end{conjecture}

\section{Integral equations}\label{finchar}

%Recall that over an algebraically closed field $F$ a complete intersection $\subset \mathbb{P}^n_F$ given by $f_1=\cdots =f_k=0$ is singular if and only if there is a point $P\in\mathbb{P}^n(F)$ such that each $f_i$ is zero at $P$ and the differentials of $f_1,\ldots, f_k$ at $P$ are linearly dependent (see e.g. Hartshorne \cite[exercise 5.8, chapter 1]{hartshorne}).

\begin{lemma}\label{definedoverz}
The hypersurface $\Sigma_{\dd,n}$ of $\Pi_{\dd,n}$ is defined by an equation with rational coefficients. Moreover, if $F$ is an arbitrary algebraically closed field, the subset $\Sigma_{\dd,n,F}\subset \Pi_{\dd,n,F}$ consisting of the equations of singular complete intersections is the zero locus of the same equation.
\end{lemma}

An element $(f_1,\ldots,f_k)$ belongs to $\Sigma_{\dd,n}$ if and only if there exist $x=(x_0,\ldots,x_n)\neq 0$
%and $(t_1,\ldots,t_k)\neq 0$
such that each $f_i$ is 0 at $x$ and the rank of the matrix $\left(\frac{\partial f_i}{\partial x_j}\right)$ at $x$ is less than $k$.
%\begin{equation}\label{conseuler}
%\sum t_i df_i|_{x}=0.
%\end{equation}
%Here $a_1,\ldots a_k$ are integers $\geq 1$ such that $a_i+\def f_i-1=\mathop{\mathrm{max}}(\deg f_1,\ldots,\deg f_k)$; they have been chosen so that each component of the equation (\label{conseuler}) is homogeneous.

%Notice that if all hypersurfaces defined by $f_i=0$ but one pass through a point $x\neq 0$ and the differentials of the $f_i$'s are linearly dependent at $x$, then the remaining hypersurface passes through $x$ as well, due to Euler's formula.
So one can obtain $\Sigma_{\dd,n}$ by taking the subvariety $\tilde\Sigma_{\dd,n}$ of $\Pi_{\dd,n}\times\C^{n+1}$ given by the equations $f_i(x)=0,i=1,\ldots,k,x\in\C^{n+1}$ and $$\mbox{a $k\times k$ minor of $\left( \frac{\partial{f}_i}{\partial x_j}\right)$ is 0}$$ %(\ref{conseuler}) with $(t_1,\ldots, t_k)\in\C^k$
and then projecting to $\Pi_{\dd,n}$. Using the classical elimination theory as presented e.g. in van der Waerden \cite[vol. II, \S 80]{waerden}, one concludes that $\Sigma_{\dd,n}$ is defined over $\Q$. More specifically, we apply theorem on pp. 158-159 there to the case when the polynomials $F_1,\ldots, F_r$ are $f_1,\ldots, f_k$ and
$k\times k$ minors of the matrix $\left( \frac{\partial{f}_i}{\partial x_j}\right)$.
The result will be a finite collection $R_1,\ldots, R_m:\Pi_{\dd,n}\to\C$ of polynomials with integer coefficients; their zero locus will be $\Sigma_{\dd,n}$. Notice that for any algebraically closed field $F$ the subvariety $\Sigma_{\dd,n,F}$ will be the zero locus  of $R_1,\ldots, R_m$ reduced modulo $F$.%We can assume in the sequel that the greatest common divisor of the coefficients of each $R_i$ is 1.

As already mentioned above, the subvariety $\Sigma_{\dd,n}$ of $\Pi_{\dd,n}$ is an irreducible hypersurface. So $R_1,\ldots, R_m$ (regarded as polynomials over $\C$) will all be multiples of a single polynomial $\triangle_{\dd,n}$, which we will assume irreducible. Moreover, e.g. by the Nullstellensatz this polynomial divides each $R_i$ and so some power of $\triangle_{\dd,n}$ is the greatest common divisor of $R_1,\ldots, R_m$. So  we may assume that this power of $\triangle_{\dd,n}$, and hence $\triangle_{\dd,n}$ itself has rational coefficients, since the coefficients of $R_1,\ldots, R_m$ are rational.
%cox little o'shea, ideals..

%This implies that $\triangle_{\dd,n}$ itself will have rational coefficients.
Multiplying if necessary $\triangle_{\dd,n}$ by a non-zero integer we can assume that all coefficients of $\triangle_{\dd,n}$ are integer and their greatest common divisor is 1. Since  each $R_i,i=1,\ldots,m$ is divisible by $\triangle_{\dd,n}$ as a polynomial over $\C$, it is also divisible by $\triangle_{\dd,n}$ as a polynomial over $\Q$, so that we can write $R_i=\frac{a}{b}f\triangle_{\dd,n}$ where $a,b$ are coprime integers with $b\neq 0$ and $f$ is a polynomial over $\Z$ with the greatest common divisor of the coefficients equal 1. By the unique factorisation this implies $b=\pm 1$, so $R_i$ is divisible by $\triangle_{\dd,n}$ as a polynomial over $\Z$.
%since finitely generated polynomial rings over $\Z$ are unique factorisation domains

This implies that for any algebraically closed field $F$ the zero locus of $\triangle_{\dd,n}$ is included in the zero locus of $R_1,\ldots, R_m$, which is $\Sigma_{\dd,n,F}$. On the other hand, $\Sigma_{\dd,n,F}$ is irreducible since it is obtained by projecting an irreducible subvariety of $\Pi_{\dd,n,F}\times F^{n+1}$ (which is defined similarly to $\tilde\Sigma_{\dd,n}$ above) to $\Pi_{\dd,n,F}$. Since $\Sigma_{\dd,n,F}$ does not coincide with the whole of $\Pi_{\dd,n,F}$, it must be equal the zero locus of $\triangle_{\dd,n}$. The lemma is proved.
$\clubsuit$

(Notice that it is possible that over some $F$ and for some $\dd$ and $n$ the equation $\triangle_{\dd,n}$ becomes reducible, and hence a power of a polynomial. This happens e.g. is when $n=1,\dd=(2),\mathop{\mathrm{char}} F=2$.)

\bigskip

In fact, generalised resultants and determinants which are described in Gelfand, Kapranov and Zelevinsky \cite{gkz}, and which go back to Arthur Cayley, give another, perhaps more explicit way to calculate $\triangle_{\dd,n}$. For general $\dd$ none of the formulae from \cite{gkz} seem to apply directly, but the general principles certainly do. Let us sketch this briefly.

Let $X$ be the quotient of $\C^{n+1}\setminus\{0\}\times \C^k\setminus\{0\}$ by the action of $\C^*\times \C^*$ given by $$(\lambda,\mu)\cdot(x,t_1,\ldots, t_k)=(\lambda x,\lambda^{a_1}\mu t_1,\ldots,\lambda^{a_k}\mu t_k)$$ where $\lambda,\mu\in\C^*$, $x\in\C^{n+1}\setminus\{0\}$ and $a_1,\ldots, a_k$ are integers $\geq 1$ such that $a_i+d_i-1=\mathop{\mathrm{max}}(d_1,\ldots,d_k)-1$. This is the total space of the projectivisation of the bundle $\sum\mathcal{O}(a_i)$ over $\mathbb{P}^n$.
%It is easy to check that $X$ is a smooth compact complex algebraic variety equipped with a map to $\mathbb{P}^n$; denote this map $p$. (Smoothness follows from the fact that the action is free and compactness from the fact that every orbit contains a point $(x,t)\in \C^{n+1}\setminus\{0\}\times \C^k\setminus\{0\}$ with $|x|=|t|=1$.)
%smooth: action is free, compact: it suffices to consider the orbits of (x,t) such that |x|=|t|=1
%this is the sum of line bundles which have sections; then look at the zero loci of the sections of those line bundles
%identify {(x,t)|x\neq 0,t\neq 0} with the total space of the bundle O(-1) pulled back to C^{n-1} as follows: (x,t)\mapsto (x,tx); the action \lambda\cdot (x,t)=(\lambda x,t\lambda^{-1}x) becomes \lambda \cdot (x,t)=(\lambda x,t)

Let $\mathcal{L}_i,i=2,\ldots,k$ be the pullback of $\mathcal{O}(d_i)$ to $X$. Moreover, the above action of $\C^*\times\C^*$ on $\C^{n+1}\setminus\{0\}\times \C^k\setminus\{0\}$ is covered by the action on $\C^{n+1}\setminus\{0\}\times \C^k\setminus\{0\}\times\C$ where $(\lambda,\mu)\in\C^*\times\C^*$ acts on the last component as multiplication by $\lambda^{\mathop{\mathrm{max}}(d_1,\ldots,d_k)-1} \mu$; the quotient $$(\C^{n+1}\setminus\{0\}\times \C^k\setminus\{0\}\times\C)/(\C^*\times \C^*)$$ is the total space of a line bundle $\mathcal{L}$ on $X$. This is $\mathcal{O}_X(1)\otimes p^*\mathcal{O}_{\mathbb{P}^n}(\mathop{\mathrm{max}}(d_1,\ldots,d_k)-1)$ where $\mathcal{O}_X(1)$ is the line bundle on $X$ constructed as in e.g. \cite[p. 160]{hartshorne}.

Set $E=\bigoplus_{i\geq 2}\mathcal{L}_i\oplus\bigoplus_{j=1}^{n+1}\mathcal{L}^{(j)}$ where $\mathcal{L}^{(j)}$'s are isomorphic copies of $\mathcal{L}$. Given an element $f=(f_1,\ldots,f_k)\in\Pi_{\dd,n}$ we obtain a section $s_f$ of $E$ by taking $f_2,\ldots,f_k$ and the $n+1$ components of $\sum t_i df_i$. Moreover, this section has a zero if and only if $f\in\Sigma_{\dd,n}$. As explained in \cite[p. 51-52]{gkz}, starting from a section $s$ of $E$ we construct the positive Koszul complex $\mathcal{K}_{+}(E,s)$.

This is a complex of sheaves on $X$; each of its terms is a sum of products of some $\mathcal{L}_i$'s and some $\mathcal{L}^{(j)}$'s. Let $\mathcal{F}$ be any such product and let us show that the higher cohomology of $\mathcal{F}$ vanishes. We have $\mathcal{F}\cong \mathcal{O}_X(b)\otimes p^*\mathcal{O}_{\mathbb{P}^n}(a)$ for some integers $a,b>0$; the (derived) pushforward $Rp_*\mathcal{O}_X(b)$ is isomorphic to $\mathop{\mathrm{Sym}}^b(\sum\mathcal{O}(a_i))$ (cf. Hartshorne \cite[exercise 8.4, chapter III]{hartshorne}). By the projection formula $Rp_*\mathcal{F}\cong \sum\mathcal{O}(a_i+a)$, which implies $H^{>0}(\mathbb{P}^n,Rp_*\mathcal{F})=H^{>0}(X,\mathcal{F})=0$.

So, by the generalised de Rham theorem, for each $s\in\Gamma E$ the exactness of $\mathcal{K}_{+}(E,s)$ is implies the exactness of its complex of the global sections $\Gamma\mathcal{K}_{+}(E,s)$, cf. \cite[Lemma 2.4]{gkz}. Let us pick a basis in $H^0$ of each term of the positive Koszul complex $\mathcal{K}_{+}(E,s)$. Given an element $f\in\Pi_{\dd,n}$, take the determinant (see e.g. \cite[Appendix A]{gkz}) of the complex of global sections of $\Gamma\mathcal{K}_{+}(E,s_f)$. The resulting map $\Pi_{\dd,n}\to\C$ will be given by a rational function with coefficients in $\Q$ that takes non-zero values on $\Pi_{\dd,n}\setminus\Sigma_{\dd,n}$. Since we already know that $\Sigma_{\dd,n}$ is irreducible, this function will be in fact polynomial defined over $\Q$ whose zero locus is precisely $\Sigma_{\dd,n}$, or the inverse of such a polynomial. Up to multiplication by an element of $\Q^*$ this polynomial is equal $\triangle_{\dd,n}$.

\bigskip

\section{Comparison with the complex case}
Set $\Sigma_{\dd,n,\Z}$ to be the scheme over $\Z$ defined by $\triangle_{\dd,n}=0$.

In order to give the proofs of theorems \ref{upperboundvectorcharp} and \ref{upboundprojectivecharp} we need to compare the \'etale cohomology with the classical one. We take an algebraically closed field $K$ of characteristic $p>0$ and a finite abelian group $A$ of order prime to $p$; moreover we assume $A$ to be a ring, so that the \'etale cohomology is equipped with an $A$-algebra structure. The variety $\Sigma_{\dd,n,K}$ introduced in the introduction is %the set of the closed points of
$\Sigma_{\dd,n,\Z}\times_{\mathop{\mathrm{Spec}}\Z} \mathop{\mathrm{Spec}}K$.

Let $R$ be the result of applying the Witt vector procedure to $K$ (see e.g. Mumford \cite[Lecture 26,
%http://eom.springer.de/W/w098100.htm -- witt vectors are complete disc valuation ring of char 0 with residue field K
%http://eom.springer.de/H/h046940.htm -- complete local rings are henselian
\S2]{mumford}); this is a complete discrete valuation ring of characteristic 0 with residue field $K$. Let $L$ be an algebraic closure of the fraction field of $R$.

Suppose we have a morphism $f:X\to Y$ of smooth schemes over $R$. Let $X_0, Y_0$, respectively $X_1,Y_1$ be the fibres of $X$ and $Y$ over $K$, respectively over $L$. The structure morphisms $X,Y\to \mathop{\mathrm{Spec}} R$ are smooth so locally acyclic %smooth morphism http://eom.springer.de/S/s085880.htm
(\cite[theorem 2.1 on p. 58]{sga45}) so we can apply the procedure explained in ibid., p. 55-56 to define the cospecialisation maps $cosp:H^{*}(X_1,A)\to H^*(X_0,A), H^{*}(Y_1,A)\to H^*(Y_0,A)$. Let us briefly recall the construction.
Since $R$ is a complete local ring,
%http://eom.springer.de/W/w098100.htm complete disc valuation
%http://en.wikipedia.org/wiki/Discrete_valuation_ring disc valuation implies local
the strict localisation of $R$ at the maximal ideal is $R$ itself.
%the elements of the complement of the max ideal of a local ring are invertible; so the map to the localisation is %iso; the map of a complete ring into the completion is iso by def, eisenbud, p 182; a completion is henselian, eisenbud p 183; when we localise and then henselise, we get the etale local ring -- the strict localisation. freitag kiehl, p. 26. If we take S=Spec R, s=X_0, then $\tilde{S}^s$ on p. 54 is just R and X_1 is L
We have the cartesian diagram (ibid, p. 56) 
$$
\begin{CD}
X_1 @>\varepsilon'>> X @<<< X_0\\
@VVV @VVV @VVV\\
\mathop{\mathrm{Spec}} L @>>> \mathop{\mathrm{Spec}}R @<<< \mathop{\mathrm{Spec}} K
\end{CD}
$$
The map $H^*(X,A)\to H^*(X_1,A)$ is an isomorphism and the cospecialisation map $H^{*}(X_1,A)\to H^*(X_0,A)$ is the composition of the inverse isomorphism and the restriction to $X_0$ (ibid, 1.6.1); similarly for the map $H^{*}(Y_1,A)\to H^*(Y_0,A)$.

It follows from the construction that the cospecialisation maps fit into the commutative diagram

\begin{equation}\label{cosp1}
\begin{CD}
 H^{*}(Y_1,A) @>cosp>> H^*(Y_0,A)\\
@VVV @VVV\\
H^{*}(X_1,A) @>cosp>> H^*(X_0,A)
\end{CD}
\end{equation}
where the vertical arrows are induced by $f$.

To prove theorems \ref{lerayhirschcharp}, \ref{upperboundvectorcharp} and \ref{upboundprojectivecharp} it would suffice to consider the case when the cardinality of $K$ is at most continuum. For theorems \ref{lerayhirschcharp}, \ref{upperboundvectorcharp} this is clear and for theorem \ref{upboundprojectivecharp} this follows from the fact that the \'etale cohomology is invariant with respect to change of base field (provided the characteristic of the coefficients is coprime with the characteristic of the base field); see \cite[Corollaire 3.3, p. 63]{sga45} or Milne \cite[Corollary 4.3, Chapter VI]{milne}.

So in the rest of the section we assume that the cardinality of $K$ is at most continuum. Then so will be the cardinality of $L$ and we can assume $L$ a subfield of $\C$. Associated to a scheme $Z$ over $\C$ are three sites: the \'etale site $Z_{\mathrm{\acute{e}t}}$, the analytic site $Z_{\mathrm{an}}$ and the topological site $Z_{\mathrm{top}}$; the covering families are formed by the \'etale maps $Z'\to Z$, local analytic isomorphisms $\tilde U\to Z(\C)$ and open embeddings $U\subset Z(\C)$ respectively. (Here $Z(\C)$ is the set of $\C$-rational points of $Z$ equipped with the structure of a complex analytic variety.)

There are natural morphisms of sites $Z_{\mathrm{\acute{e}t}}\gets Z_{\mathrm{an}}\to Z_{\mathrm{top}}$. M. Artin's comparison theorem (\cite[expos\'e XI]{sga43}) states that the induced cohomology maps with constant finite coefficients are isomorphisms. Moreover, since these isomorphisms are induced by maps of sites, they commute with the cohomology maps induced by morphisms of schemes over $\C$.

Applying this to $X_1, Y_1$ and $f$ as above we can extend diagram (\ref{cosp1}) to
\begin{equation}\label{compfunc}
\begin{CD}
H^*(Y_1(\C),A) @>{\cong}>> H^{*}(Y_1,A) @>cosp>> H^*(Y_0,A)\\
@VVV @VVV @VVV\\
H^*(X_1(\C),A) @>{\cong}>> H^{*}(X_1,A) @>cosp>> H^*(X_0,A)
\end{CD}
\end{equation}
where $X_1(\C)$ and $Y_1(\C)$ are the complex analytic varieties obtained by taking the sets of $\C$-rational points of $X_1$ and $Y_1$ (base changed to $\C$) respectively and the vertical arrow on the left is the cohomology map induced by the continuous map of topological spaces (obtained from $f$). The horizontal compositions in this diagram will be denoted $comp_X$ and $comp_Y$ respectively and will be called the {\it comparison morphisms}.

These are not isomorphisms in general; the next step in the proof of theorems \ref{lerayhirschcharp}, \ref{upperboundvectorcharp} and \ref{upboundprojectivecharp} is to show that they are when the source is the cohomology of $\C^m\setminus\{0\}$, $\SL_m(\C)$ or $\GL_m(\C)$ with coefficients in $A$. More precisely, to construct these maps we first need to define the corresponding smooth schemes over $R$. This can be done as follows. Consider the complement of the origin in $\mathbb{A}^m_R$ (i.e. the complement of the closed subscheme defined by $x_1=\cdots =x_m=0$), the closed subscheme $\SL_{m,R}$ of $\mathbb{A}^{m^2}_R$ (the affine space of $m\times m$ matrices over $R$) defined by the equation $\mathop{\mathrm{det}}=1$ and complement $\GL_{m,R}$ in $\mathbb{A}^{m^2}_R$ of the closed subscheme defined by the equation $\mathop{\mathrm{det}}=0$.

From now on we assume that $A$ is the ring $\mathbb{Z}/r$ with $r$ coprime with $p$.

\begin{lemma}\label{cohglmcharp}
The resulting comparison morphisms $$H^*(\C^m\setminus\{0\},A)\to H^*(K^m\setminus\{0\},A),$$ $$H^*(\SL_m(\C),A)\to H^*(\SL_m(K),A),$$ $$H^*(\GL_m(\C),A)\to H^*(\GL_m(K),A)$$ are isomorphisms.
\end{lemma}

The proof is by computing the left and the right hand sides simultaneously, and noticing that both in the complex case and in the \'etale case the computations give identical results. This is probably standard but we will give a sketch however, since the details seem not to be available in the literature.

We start with $\C^n\setminus\{0\}$. %The scheme over $R$ that we need to consider is $\mathbb{A}^m_R\setminus\{0\}$, the affine space minus the origin.
Let $X$ be the complement of $\mathbb{A}^m_R\setminus\{0\}$ in $\mathbb{A}^m_R$ and let $i$ be the closed embedding $X\subset\mathbb{A}^m_R$.  By the purity theorem (see e.g. Milne \cite[Theorem VI 5.2]{milne} or \cite[p. 63]{sga45}), we have $i^!\underline{A}_{\mathbb{A}^m_R}=i_*\underline{A}_{X}(-m)[-2m+1]$ where $\underline{A}_Y$ is the constant sheaf with stalk $A$ on a scheme $Y$, the round brackets stand for the Tate twist and the square brackets indicate a degree shift. Since $A=\Z/r$ with $r$ non-divisible by $p$ Tate twisting the constant sheaf gives the same sheaf, up to a non-canonical isomorphism.
%milne p. 125
The cospecialisation map is an isomorphism for $\underline{A}_{\mathbb{A}^m_R}$ and $i^!\underline{A}_{\mathbb{A}^m_R}$, so it is an isomorphism for $Rj_*\underline{A}_{\mathbb{A}^m_R}$ where $j:\mathbb{A}^m_R\setminus\{0\}\to \mathbb{A}^m_R$ is the open embedding. This proves the statement of the lemma for the affine space minus the origin.

One then proceeds as follows. Let us consider the case of $\GL_m$. Let $\GL_{m,K}$ and $\GL_{m,L}$ be the closed and the generic fibers of $\GL_{m,R}$ (i.e. the fibers over $\mathop{\mathrm{Spec}} K$ and $\mathop{\mathrm{Spec}} L$ respectively). Since $H^*(\GL_{m,L},A)\cong H^*(\GL_{m}(\C),A)$ by M. Artin's comparison theorem, we only need to see what happens with the cospecialisation map $H^*(\GL_{m,L},A)\to H^*(\GL_{m,K},A)$.

One computes the cohomology of $\GL_{m,K}$ and $\GL_{m,L}$ separately; as in section 2 of \cite{division} one finds by induction a system of canonical multiplicative generators in odd degrees such that all generators but the last one restrict to the canonical generators of $\mathrm{GL}_{m-1}$ and the last one is the pullback of a generator of the "sphere" $\mathbb{A}^m$ minus the origin. The generators are defined modulo multiplication by units of $A$. Then one proves by induction using the functoriality (with respect to the mappings $\mathrm{GL}_{m-1}\to\mathrm{GL}_m$ and $\mathrm{GL}_m\to\mathbb{A}^m$ minus the origin; both these mappings are defined over $R$) that the comparison maps are isomorphisms.

The case of $\mathrm{SL}_m$ is similar.$\clubsuit$

A different proof of the same result is sketched in \cite[p. 230]{sga45}.
%matsumura corollary 2 p 152; eisenbud-harris p. 37 atiyah-macdonald exercise 2.15

Let $\cc_1^m(K,A),\ldots,\cc_m^m(K,A)$ be the classes in $H^*(\GL_{m,K},A)$ that correspond to the reduction modulo $A$ of the classes $\cc_1^m,\ldots,\cc_m^m$ from section  2 in \cite{division} via the isomorphism of lemma \ref{cohglmcharp}.
(Here the reduction modulo $A$ is the map from the integral cohomology to the cohomology with $A$ coefficients induced by the unit $\Z\to A$.)\label{redmoda} Let $\cc_i^m(K,\Q_\ell)$ be the classes obtained from $\cc_i^m(K,\Z/\ell^r)$ by taking the limit as $r\to\infty$ and tensoring with $\Q_\ell$.

Our next task in this section is to construct the orbit maps over $R$.
Take a nonsingular multidegree $\dd$ complete intersection in $\mathbb{P}^n_K$ given by $f_1=\cdots=f_k=0$ and let $f$ be a lift of $(f_1,\ldots,f_k)$ to $R$. There is the orbit map $\mathbb{A}^{n^2}_R\to\Pi_{(d),n,R}$ obtained by ``transforming $f$ by linear substitutions''. We would like to compute the pullback $\tilde\triangle$ of $\triangle_{\dd,n}$ under this map. Recall that we consider $R$ as a subring of $\C$; over the complex numbers the zero locus of $\tilde\triangle$ is the same as the zero locus of the determinant hypersurface, so by the Nullstellensatz we have $\tilde\triangle=a\mathrm{det}^l$ where $a$ is a non-zero complex number and $l$ is an integer $\geq 1$. By evaluating at the identity matrix we see that $a=\triangle_{\dd,n}(f)$. Since $R$ is a local ring and $f$ is obtained by lifting the equations of a non-singular complete intersection, we see that $a\in R$ and moreover that this is a unit of $R$.

So we get the orbit map from the complement of the $\mathrm{det}=0$ subscheme of $\mathbb{A}^{n^2}_R$ to the complement of the $\triangle_{\dd,n}=0$ subscheme of $\Pi_{\dd,n,R}$, i.e. a map

\begin{equation}\label{orbitcharp}
\GL_{n+1,R}\to \Pi_{\dd,n,R}\setminus\Sigma_{\dd,n,R}.
\end{equation}

Moreover, since the origin in $\Pi_{\dd,n,R}$ is a subscheme of $\Sigma_{\dd,n,R}$ we can compose the above map with the natural morphism from $\Pi_{\dd,n,R}$ minus the origin to the projectivisation $\mathbb{P} \Pi_{\dd,n,R}$ of $\Pi_{\dd,n,R}$.
%on the complement of the origin the coordinate functions generate the stalks of $\mathcal{O}$ as modules over themselves so by Hartshorne theorem 7.1 p. 150 we get a map to the projective space
Moreover, we get in fact a morphism from $\GL_{n+1,R}$ to the complement of the projectivisation $\mathbb{P}\Sigma_{\dd,n,R}$ of $\Sigma_{\dd,n,R}$ in $\mathbb{P} \Pi_{\dd,n,R}$ i.e. the complement of the closed subscheme defined by the homogeneous ideal generated by $\triangle_{\dd,n}=0$.

Recall from \cite{division}, section 2, that there are classes $\aaa^{\dd,n}_i\in \in H^{2i-1}(\Pi_{\dd,n}\setminus\Sigma_{(d),n},\Z)$ which, when pulled back under any orbit map $\GL_{n+1}(\C)\to \Pi_{\dd,n}\setminus\Sigma_{(d),n}$ give $$\m_i^{\dd,n}\cc_i^{n+1}$$.

Now take $Y=\Pi_{\dd,n,R}\setminus\Sigma_{\dd,n,R}, X=\GL_{n+1,R}$ and apply the comparison map $comp_Y$ to the classes $\aaa^{\dd,n}_i$  reduced modulo $A$ (as above on page \pageref{redmoda}, reduction modulo $A$ is the cohomology map induced by the unit $\Z\to A$). Denote the resulting classes $\aaa^{\dd,n}_1(K,A),\ldots,\aaa^{\dd,n}_{n+1}(K,A)$; they live in the cohomology of $\Pi_{\dd,n,K}\setminus\Sigma_{\dd,n,K}.$ Using the commutativity of diagram (\ref{compfunc}) and the explicit formulae for $\m_i^{\dd,n}$ from \cite[sections 3 and 4]{division} we see that the pullback of each $\aaa^{\dd,n}_i(K,A)$ under the orbit map (\ref{orbitcharp}) is $\m_i^{\dd,n}\cc_i^{n+1}(K,A)$. By taking the limit of $\aaa^{\dd,n}_1(K,\Z/\ell^r)$ as $r\to\infty$ and tensoring with $\Q_\ell$ we get classes in $H^*(\Pi_{\dd,n,R}\setminus\Sigma_{\dd,n,R},\Q_\ell)$ whose pullbacks under the orbit maps are $\m_i^{\dd,n}\cc_i^{n+1}(K,\Q_\ell)$.

Now we turn to the case when $X$ is a hypersurface of degree $d>2$. We can apply the above procedure to $Y$ equal $\mathbb{P}\Pi_{(d),n,R}\setminus \mathbb{P}\Sigma_{(d),n,R}$ and $X=\SL_{n+1,R}$. We take $comp_Y$ of the reduction modulo $A$ of the classes $\aaa_i^{(d),n,proj},i=2,\ldots,n+1$ introduced at the end of section 5 in \cite{division} to get classes $\aaa_i^{(d),n,proj}(K,A)$ that live in $H^*(\mathbb{P}\Pi_{(d),n,K}\setminus \mathbb{P}\Sigma_{(d),n,K},A)$. Again using the commutativity of (\ref{compfunc}) and the formula in the complex case at the end of section 5, \cite{division} we conclude that $\aaa_i^{(d),n,proj}(K,A)$ pulls back to
\begin{equation}\label{projexplformcharp}
\frac{\LCM(C^{n-i+1}_{n+1}(d-1)^{n-i+1},(n+1)(d-1)^n)}{C^{n-i+1}_{n+1}(d-1)^{n-i+1}}\m_i^{(d),n}\cc_i^{n+1}(K,A)
\end{equation}
under the orbit map $\SL_{n+1,R}\to \mathbb{P}\Pi_{(d),n,R}\setminus \mathbb{P}\Sigma_{(d),n,R}$.

Let us state the results of the last two paragraphs as a lemma, since we will need them later.

\begin{lemma}\label{explcoeffscharp}
1. For $i=1,\ldots, n+1$ there exist classes $\aaa^{\dd,n}_i(K,A)\in H^{2i-1}(\Pi_{\dd,n,K}\setminus\Sigma_{(d),n,K},A)$ and $\aaa^{\dd,n}_i(K,\Q_\ell)\in H^{2i-1}(\Pi_{\dd,n,K}\setminus\Sigma_{(d),n,K},\Q_\ell)$ that, for some choice of a point of $\Pi_{(d),n,K}\setminus\Sigma_{(d),n,K}$, pull back to $\m_i^{\dd,n}\cc_i^{n+1}(K,A)$, respectively $\m_i^{\dd,n}\cc_i^{n+1}(K,\Q_\ell)$, under the resulting orbit map $\GL_{n+1}(K)\to \Pi_{\dd,n,K}\setminus\Sigma_{\dd,n,K}$.

2. For $i=2,\ldots, n+1$ there exist classes $\aaa^{(d),n,proj}_i(K,A)\in H^{2i-1}(\mathbb{P}\Pi_{(d),n,K}\setminus\mathbb{P}\Sigma_{(d),n,K},A)$ that, for some choice of a point of $\mathbb{P}\Pi_{(d),n,K}\setminus\mathbb{P}\Sigma_{(d),n,K}$, pull back to (\ref{projexplformcharp}) under the resulting orbit map $\SL_{n+1}(K)\to\mathbb{P}\Pi_{(d),n,K}\setminus\mathbb{P}\Sigma_{(d),n,K}$.
\end{lemma}

(We will see shortly that the choice of the points for the orbit maps in cohomology does not matter, just as it doesn't in the complex case.)

$\clubsuit$

In a similar way one can take $Y=\GL_{n+1,R}$; then $X=\GL_{1,R}$ acts on $Y$ by multiplication by scalar matrices. Applying $comp_Y$ to the reduction of $\cc_1^{n+1}\in H^1(\GL_{n+1}(\C),\Z)$ modulo $A$ we obtain

\begin{lemma}\label{explcoeffscharpgroup}
There exists a class $\cc_1^{n+1}(K,A)\in H^1(\GL_{n+1}(K),A)$ that pulls back to a generator of $H^1(K^*,A)$ under the orbit map for some choice of a point in $\GL_{n+1}(K)$.
\end{lemma}

$\clubsuit$

\section{Stabilisers}

Here we prove some statements on automorphisms and vector fields on complete intersections. We start with vector fields.

\begin{lemma}\label{novectorfields}
A smooth complete intersection $X$ of multidegree $\dd=(d_1,\ldots,d_k)$ in $\mathbb{P}^n_F$ where $F$ is an algebraically closed field admits no non-zero tangent vector fields, unless $X$ is a quadratic hypersurface, a cubic curve in $\mathbb{P}^2_F$, the intersection of two quadrics in $\mathbb{P}^3_F$, or the intersection of two quadrics in $\mathbb{P}^n_F, n\geq 5$ odd and $char(F)=2$.
\end{lemma}

This was proved by O. Benoist in \cite[Th\'eor\`eme 3.1]{benoist}.

\begin{lemma}\label{cornovectorfields}
There are no vector fields $\in H^0(\mathbb{P}^n_F,T\mathbb{P}^n_F)$ that are tangent to a smooth complete intersection $X$ of multidegree $\dd=(d_1,\ldots,d_k)$, except when $X$ is a quadratic hypersurface, a cubic curve in $\mathbb{P}^2_F$ and $char(F)=3$, or the intersection of two quadrics in $\mathbb{P}^n_F, n\geq 3$ odd and $char(F)=2$.
\end{lemma}

{\bf Proof.} %Indeed, if such a vector field existed, then from lemma \ref{novectorfields} it would follow that the field would be zero when restricted to $X$. The zeroes of a vector field $\in H^0(\mathbb{P}^n_F,T\mathbb{P}^n_F)$ form a union of disjoint projective subspaces. When $\dim X>0$ the corollary follows from the fact that $X$ is connected (see e.g. Hartshorne \cite[chapter III, excercise 5.5]{hartshorne}) and does not lie in a hyperplane.
The case $\dim X>0$ is covered by theorem 3.1 in \cite{benoist}. It remains to consider the case when $\dim X=0$.
Let $X$ be a complete intersection in $\mathbb{P}^n_F$ given as $f_1=\cdots=f_k=0$ with $(f_1,\ldots, f_k)\in\Pi_{\dd,n,F}\setminus\Sigma_{\dd,n,F}$ and $k=n$. We have to show that there is no non-zero tangent vector field $\in H^0(\mathbb{P}^n_F,T\mathbb{P}^n_F)$ whose zero locus contains $X$. The case when $n=1$ is clear, so in the sequel we assume $k=n>1$.

The zeroes of a non-zero element of $H^0(\mathbb{P}^n_F,T \mathbb{P}^n_F)$ form a family of projective subspaces of $\mathbb{P}^n_F$ that can be included in a union $L_1\sqcup L_2$ of non-intersecting projective subspaces. Let $Y\supset X$ be the subvariety of $\mathbb{P}^n_F$ (in fact, a curve) given as the zero locus of $f_1,\ldots,f_{k-1}$; $Y$ can not be contained in a hyperplane (or else $X$ would be as well, which would contradict the assumption that all $d_i\geq 2$. So, since $Y$ is connected, it can not be contained in $L_1\sqcup L_2$. Let $P$ be a point on $Y$ outside $L_1\sqcup L_2$ and let $H_i$ be a hyperplane through $P$ and $L_i$. By B\'ezout's theorem, $Y$ intersects each $H_i$ in $\leq \Pi_{i=1}^{k-1} d_i$ points; moreover, at least one of the points of $Y\cup H_1$, namely $P$, also belongs to $Y\cup H_2$. So $Y\cup (H_1\cap H_2)$ contains $<2\Pi_{i=1}^{k-1} d_i\leq \Pi_{i=1}^{k} d_i$ points. Since $X\subset Y\cap (L_1\cup L_2)\subset Y\cap (H_1\cup H_2)$, we conclude that $X$ contains less than $\Pi_{i=1}^{k} d_i$ points, which contradicts the assumption that $X$ is non-singular.$\clubsuit$
\begin{Prop}\label{finitestab}
Let $F$ be an algebraically closed field and let $\dd\neq (2)$. Then the stabilisers both of an element of $(f_1,\ldots,f_k)\in\Pi_{\dd,n,F}\setminus\Sigma_{\dd,n,F}$ and of the corresponding multidegree $\dd$ complete intersection in $\mathbb{P}^n_F$ are finite.
\end{Prop}

Note that here we consider stabiliser subgroups; these are subgroups of $\GL_{n+1}(F)$, respectively $\PGL_{n+1}(F)$.
%the sets of the closed points of the stabiliser subschemes.

{\bf Proof.} The case $k<n$ has been treated by O. Benoist in \cite[Th\'eor\`eme 3.1]{benoist}. (Whenever we can apply lemma \ref{novectorfields} we do so and the exceptional cases have to be considered separately.) It remains to prove the proposition when $k=n$ or $k=n+1$. In the first case the proposition follows from lemma \ref{cornovectorfields}.

Suppose $k=n+1$. Since $(f_1,\ldots,f_{n+1})\in\Pi_{\dd,n,F}\setminus\Sigma_{\dd,n,F}$, the zero locus of the ideal generated by $f_1,\ldots, f_{n+1}$ is the origin. Let $r_i$ be $\frac{1}{d_i}\LCM(d_1,\ldots,d_{n+1})$. The morphism $f:\mathbb{A}^{n+1}_F\to\mathbb{A}^{n+1}_F$ with components $f^{r_1}_1,\ldots, f^{r_{n+1}}_{n+1}$ is finite and surjective. (Indeed, it is the restriction to $\mathbb{A}^{n+1}_F=\bar f^{-1}\mathbb{A}^{n+1}_F$ of a non-constant (hence finite and surjective) morphism $\bar f:\mathbb{P}^{n+1}_F\to\mathbb{P}^{n+1}_F$, so $f$ is also finite, hence closed and hence surjective since it is dominant.)
%milne p. 4
%So the morphism $f:\mathbb{A}^{n+1}_F\to\mathbb{A}^{n+1}_F$ with components $f_1,\ldots, f_{n+1}$ is also finite, hence closed,
%hartshorne p. 91, p.105 proper
%and hence surjective (since it is dominant).

Now we can construct $n+1$ linearly independent points $x_0,\ldots,x_n$ in $\mathbb{A}^{n+1}_F$ such that the images $f(x_0),\ldots,f(x_n)$ are pairwise distinct. (This can be done as follows: take any non-zero point for $x_0$, then take $x_1\not\in$ the preimage of the image of the line spanned by $x_0$, then take $x_2\not\in$ the preimage of the image of the 2-plane spanned by $x_0$ and $x_1$ and so on.) The stabiliser $\in\GL_{n+1}(F)$ of $(f_1,\ldots,f_{n+1})$ preserves each finite set $f^{-1}(f(x_i)),i=0,\ldots,n$, so it has a finite-index subgroup that preserves each $x_i$. This subgroup is the identity, so the stabiliser is finite. The proposition in proved.$\clubsuit$

It may be possible to show directly that the stabilisers of the elements of $\Pi_{\dd,n}\setminus\Sigma_{\dd,n}$ do not contain unipotent transformations, and then to deduce proposition \ref{finitestab} from theorem \ref{main} and its \'etale analogue, theorem \ref{lerayhirschcharp} without using the non-existence of algebraic vector fields on smooth complete intersections. However, we will not attempt this in this paper.

\bigskip

Let us mention here one more related result that we will need in the sequel. The quotient stacks $[(\Pi_{\dd,n}\setminus\Sigma_{\dd,n})/\GL_{n+1}(F)]$ are in general not Deligne-Mumford when $char(F)>0$, which poses some problems. However, there are separated Deligne-Mumford stacks that will be just as good for our purposes.

Starting from $\Pi_{\dd,n}\setminus\Sigma_{\dd,n}$ we can construct an open subvariety $U$ of $\Pi_{i=1}^k \mathbb{P}^{D_i-1}$ where $D_i$ is the dimension of the space of homogeneous polynomials of degree $d_i$ is $n+1$ variables. The group $\PGL_{n+1}(F)$ acts on $U$ in a natural way.

\begin{lemma}
The quotient stack $[U/\PGL_{n+1}(F)]$ is Deligne-Mumford and separated, except when $\dd=(2)$ or $\dd=(3), char(F)=3,n=2$, or $\dd=(2,2),char(F)=2, n$ is odd.
\end{lemma}

When $\dd=(2)$ or $\dd=(3), char(F)=3,n=2$ it is impossible for $[U/\PGL_{n+1}(F)]$ to be Deligne Mumford as in these cases complete intersections have non-zero tangent vector fields that come from vector fields on $\mathbb{P}^n_F$ (lemma \ref{cornovectorfields}), which means that the projective automorphism group scheme is infinite or non-reduced. Moreover, a generic complete intersection of two quadrics in an odd-dimensional projective space over a field of characteristic $2$ can be given by equations $$\sum_{i=0}^{r-1}x_i x_{i+r}=0, \sum_{i=0}^{r-1}a_i x_i x_{i+r}+\sum_{i=0}^n b_i x_i^2=0$$
where $r=\frac{n+1}{2}$ and $a_i,b_i\in F$, see \cite[proof of Proposition 3.8]{benoist}. For any choice of $c_0,\ldots, c_{r-1}$ the vector field $\sum_{i=0}^{r-1} c_i\left(x_i \frac{\partial}{\partial x_i}+x_{i+r}\frac{\partial}{\partial x_{i+r}}\right)$ is tangent to both the above quadrics, so $[U/\PGL_{n+1}(F)]$ can't be Deligne-Mumford.

In all the remaining cases the Deligne-Mumford property follows from lemma \ref{cornovectorfields}. The separatedness of $[U/\PGL_{n+1}(F)]$ is equivalent to the properness of the action of $\PGL_{n+1}(F)$ on $U$, cf. \cite[Proposition 4.17]{edidin}. To show that $\PGL_{n+1}$ acts properly on $U$ it would suffice to prove that every point of $U$ is properly stable with respect to some line $\PGL_{n+1}(F)$-bundle $L$ on $U$, \cite[Corollary 2.5]{git}.
%git, def 1.8; cf. def 1.7
But $U$ is in fact an affine variety as it is the quotient of the affine variety $\Pi_{\dd,n}\setminus\Sigma_{\dd,n}$ with respect to an action of $(F^*)^k$.
So we can take $L$ to be the trivial bundle; then every element of $U$ will be properly stable since all orbits of $\PGL_{n+1}(F)$ have the same dimension (and hence, are closed).$\clubsuit$

\section{Proofs of the theorems}

Now we need an \'etale analogue of the observation on p. \pageref{obser}. It is given by the following lemma.

\begin{lemma}\label{finstabcharp}
Let $G=\GL_r(K)$ or $\SL_r(K)$ acting with finite stabilisers on a connected affine variety $X$ over $K$. Set $m=r^2$ if $G=\GL_r(K)$ and $m=r^2-1$ if $G=\SL_r(K)$. Then all orbit maps induce the same map in \'etale cohomology with coefficients in $\mathbb{Z}/\ell^r$. Suppose that for some $x\in X$ the index of the image in $H^m(G,\Z/\ell^r)$ of the cohomology map induced by the orbit map of $x$ contains an element that spans a non-zero subgroup of index $a\in\Z$. Then the $\ell$-part of the stabiliser of each point of $X$ is divisible by $a$.
\end{lemma}

Let us first show that all orbit maps induce the same map in \'etale cohomology.\label{samemap} Since they are all compositions of ``horizontal'' inclusions $G\to G\times X$ and the action map $G\times X\to X$, it suffices to show that any two such inclusions of induce the same cohomology map. But this follows from the connectedness of $X$,
%freitag-kiehl p. 40: since $X$ is connectedthe group of sections of the constant $A$-sheaf is iso to $A$, so $H^0(X,A)=A$; $H^*(G,A)\otimes H^*(X,A)\to H^*(G\times X,A)$ is iso by K\"unneth. Each class $c\otimes 1$ where $c\in H^*(G,A)$ restricts to $c$ under any inclusion -- since the composition inclusion first, then projection is the identity. Moreover, under any inclusion each class $c\otimes c'$ where $c\in H^*(G,A)$ and $c'\in H^{>0}(X,A)$ restricts to 0 since each class $1\times c'$ with $c'\in \in H^{>0}(X,A)$ restricts to 0; the latter follows from the fact that it is induced by inclusion $G\to G\times X$ and the projection to $X$ and the fact that the cohomology of a point is trivial.
the fact that the cohomology of $G$ is a free $\Z/\ell^r$-module (as proved above) and the K\"unneth formula (see e.g. Milne \cite[Remark 8.25 on p. 267]{milne}; cf. \cite[Th\'eor\`emes de finitude en cohomologie $\ell$-adique]{sga45}).

Starting from this point the proof follows very closely the proof of the observation given in section 5 of \cite{division}: we have an inclusion of subgroups of $H^*(G,\Z/\ell^r)$:
\begin{equation}\label{3groups}
\mbox{a subgroup of index $a$}\subset\mathop{\mathrm{Im}} orb_x^*\subset \mathop{\mathrm{Im}} p^*
\end{equation}
where $p:G\to G/Stab (x)$ is the projection. However there is one important difference: we can no longer use maximal compact subgroups in the argument to show that the index of the image of $p_*$ is equal to the $\ell$-part of $Stab(x)$. Instead we use the following lemma.

\begin{lemma}\label{fingroup}
Let $Y$ be an affine scheme over $K$ of dimension $m$, $H$ be a finite group acting on $Y$ and let $p:Y\to Y/H$ be the projection. Assume that the group $H^m(Y,\Z/\ell^r)$ is isomorphic to $\Z/\ell^r$ and that the action of  $H$ on $H^m(Y,\Z/\ell^r)$ is trivial. Then

1. If $\ell^r$ is greater than the $\ell$-part of $\#H$, then the index of $p^*(H^m(Y/H,\Z/\ell^r))$ in $H^*(Y,\Z/\ell^r)$ is equal the $\ell$-part of $\#H$.

2. If $\ell^r$ is smaller than or equal the $\ell$-part of $\#H$, then $p^*(H^m(Y/H,\Z/\ell^r))=0$.
\end{lemma}
%quotients of affine eschemes by finite groups exist, see appendix, mustata

{\bf Proof\footnote{Based on a suggestion by Bas Edixhoven.}.} Set $F$ to be the constant $\Z/\ell^r$-sheaf on $Y/H$. Since the higher derived images of $p^{-1}F$, the constant $\Z/\ell^r$-sheaf on $Y$, under $p$ vanish,
%can check this on the stalks using proper base change; closed geometric points suffice -- milne p. 65; fef of a geometric point -- milne p. 60
we have a natural isomorphism $H^*(Y,\Z/\ell^r)\cong H^*(Y/H,p_*p^{-1}F)$.

We have the restriction map $rest:F\to p_*p^{-1}F$ and the trace map $tr:p_*p^{-1}F\to F$. The former of these is obtained from the pair $(p^{-1},p_*)$ of adjoint functors and gives us the pullback $H^*(Y/H,F)\to H^*(Y,\Z/\ell^r)$ using the above identification. The latter is constructed as follows. If $U\to Y/H$ is \'etale, then a section of $p_*p^{-1}F$ over an \'etale $U\to Y/H$ can be viewed as a function from the set of connected components of $Y\times_{Y/H} U$ to $\Z/\ell^r$, and a section of $F$ over $U$ as a function from the set of connected components of $U$ to $\Z/\ell^r$. We construct a map $H^0(U,p_*p^{-1}F)\to H^0(U,F)$ by taking $tr(s)(C)=\sum_{g\in H} s(g\cdot C')$ where $s\in H^0(U,p_*p^{-1}F)$, $C$ is a connected component of $U$ and $C'$ is a component in $Y\times_{Y/H} U$ that maps to $C$. From this we get a map of sheaves, which we again denote as $tr$. From this definition it follows that $tr\circ rest$ is multiplication by the order of $H$.

Notice that the trace map is surjective on the stalks so there is a sheaf $Q=\ker (tr)$ such that the sequence $$0\longrightarrow Q\longrightarrow p_*p^{-1}F\stackrel{tr}{\longrightarrow} F\longrightarrow 0$$ is exact. The quotient $Y/H$ is again affine, so the group $H^{m+1}(Y/H,Q)=0$ (see e.g. \cite[expos\'e XIV]{sga43}) and so the trace map induces a surjection $H^m(Y,\Z/\ell^r)\to H^m(Y/H,F)$.

The former of these two groups is assumed to be $\Z/\ell^r$, so the latter is $\Z/\ell^s,s\leq r$. Slightly abusing the terminology we will use ``trace'' and ``restriction'' to denote the resulting cohomology maps. We can assume that the trace map in degree $m$ is the standard map $\Z/\ell^r\to\Z/\ell^s$ i.e. that it takes the additive generator $1\in\Z/\ell^r$ to the additive generator $1\in\Z/\ell^s$. Under the restriction map $1$ goes to $\ell^c b,r\geq c\geq r-s$ where $b$ is invertible modulo $\ell^r$; when we take the trace of this we get $\ell^c b'\in\Z/\ell^s$ where $b'$ is invertible modulo $\ell^s$. Clearly, $\ell^c$ is the index of the image of the restriction map in degree $m$, or, which is the same, the index of $p^*(H^m(Y/H,\Z/\ell^r))$ in $H^m(Y,\Z/\ell^r)$.

The composition $tr\circ rest$ acts as multiplication by $\#H$ on $H^m(Y,\Z/\ell^r)$. To see this note that this map takes an element to the sum of all its $H$-images and that the action of $H$ on $H^m(Y,\Z/\ell^r)$ is assumed to be trivial. So the multiplication by $\#H$ induces the same map as the multiplication by $\ell^c b$ on $H^m(Y,\Z/\ell^r)$. This implies the lemma.$\clubsuit$
%But the composition $tr\circ adj$ is the multiplication by $\#H$ (one can check this on the stalks), so $\ell^a$ is also the $\ell$-part of $\#H$. The lemma is proved.$\clubsuit$

Now we can return to the proof of lemma \ref{finstabcharp}. We would like to apply lemma \ref{fingroup} to $Y=G$ with $H$ being the stabiliser of a point of $X$. To do so we must first check that $G$ acts identically on its cohomology. This follows from the K\"unneth formula, cf. the proof of the fact that all orbit maps induce the same cohomology map on p. \pageref{samemap}.

Now, since the smallest group in (\ref{3groups}) is non-zero by the hypothesis of the lemma, we apply lemma \ref{fingroup} to deduce that the index of the largest one in (\ref{3groups}) is equal to the $\ell$-part of $Stab(x)$. This divides the index of the smallest group, which is $a$. Lemma \ref{finstabcharp} is proved.$\clubsuit$

Theorem \ref{upperboundvectorcharp} now follows from lemma \ref{finstabcharp} and lemma \ref{explcoeffscharp} (we apply lemma \ref{finstabcharp} to the action of $\GL_{n+1}(K)$ on $\Pi_{(d),n,K}\setminus\Sigma_{(d),n,K}$ with $a$ being the product of $\m_i^{(d),n}$ for $i=1,\ldots,n+1$ and $\ell$ being a prime number different from the characteristic of $K$).

In a similar way we can deduce the statement of theorem \ref{upboundprojectivecharp} on the projective automorphism groups by applying lemma \ref{finstabcharp} to the action of $\SL_{n+1}(K)$ on $\mathbb{P}\Pi_{(d),n,K}\setminus\mathbb{P}\Sigma_{(d),n,K}$ (we use that fact that for any prime $\ell$ different from the characteristic of $K$ the $\ell$-part of the kernel of $\SL_{n+1}(K)\to\PGL_{n+1}(K)$ is equal the $\ell$-part of $n+1$). To complete the proof of \ref{upboundprojectivecharp} we need to show that for non-quadratic hypersurfaces of dimension $>1$ the projective automorphism group coincides with the whole automorphism group. This is a consequence of the following lemma.

\begin{lemma}\label{projautisaut}
Let $V$ be a smooth complete intersection of dimension $>1$ and multidegree $\dd$ in $\mathbb{P}^n_K$.
The automorphism group of $V$ coincides with the projective automorphism group except in the following cases: $d=(4),n=3$; $\dd=(2,3),n=4$; $\dd=(2,2,2),n=5$. (The exceptions are the complete intersection $K3$ surfaces.) 
\end{lemma}

The idea of the proof is to interpret the class of the hyperplane section divisor $\in \mathop{\mathrm{Pic}}V$ in an invarint way. This interpretation comes from two sources.

First, when $\dim V>2$ the group $\mathop{\mathrm{Pic}}V$ is $\Z$ (see \cite[expos\'e XII, corollaire 3.7]{sga2}) and the hyperplane section generates the effective cone. Hence the statement of the lemma for complete intersections of dimension $>2$.

Second, the canonical class of $V$ is $\sum_{i=1}^k d_i-n-1$ times the hyperplane section, see e.g. Hartshorne \cite[exercise 8.4, Chapter II]{hartshorne}. So if $\sum_{i=1}^k d_i-n-1\neq 0$ and $\mathop{\mathrm{Pic}}V$ is torsion free, we are done. The group $\mathop{\mathrm{Pic}}V$ is indeed torsion free when $V$ is a surface (see \cite[Th\'eor\`eme 1.8 and lemme 1.9, p. 49]{sga72}). This proves the lemma for non $K3$ surfaces (and gives an alternative proof for higher dimensional complete intersections with non-zero canonical class).
$\clubsuit$
%Badescu, Nagoya Math Journ, Vol 71, 1978, 169-179

Let us now give a few remarks on conjecture \ref{lerayhirschcharp}. %We start with a preliminary lemma.

There are two main difficulties compared to the complex case. The first one is that there seems to be no slice theorem that would be readily applicable to the action of $\GL_{n+1}(K)$ on $\Pi_{\dd,n}\setminus\Sigma_{\dd,n}$, so we have to proceed in a roundabout way and prove the stack theoretic version first. The second difficulty is that the quotient stacks $[(\Pi_{\dd,n}\setminus\Sigma_{\dd,n})/\GL_{n+1}(K)]$ are in general not Deligne-Mumford: e.g. when all $d_i$'s are divisible by $char(K)$, the stabiliser of each elelement of $\Pi_{\dd,n}\setminus\Sigma_{\dd,n}$ contains a non-reduced subgroup. To remedy this we consider a slightly different stack, which is Deligne-Mumford and separated, has the same cohomology as $[(\Pi_{\dd,n}\setminus\Sigma_{\dd,n})/\GL_{n+1}(K)]$, and whose coarse moduli space coincides with $(\Pi_{\dd,n}\setminus\Sigma_{\dd,n})/\GL_{n+1}(K)$.

Since we will only be interested in quotient stacks here, our main reference for stacks will be Bernstein-Lunts \cite{berlunts}. As noticed in \cite[4.3]{berlunts}, all constructions from the first three chapters of \cite{berlunts} can be extended to the algebraic case with obvious modifications, which we indicate below.

Let $X$ be a variety over $K$ and let $G$ be an algebraic group acting on $X$. We assume that $G$ is a closed and reduced subgroup of some $\GL_{m}(K)$. Recall that the quotient stack $[X/G]$ is the category whose objects are $G$-torsors $P\to T$ equipped with $G$-equivariant maps $f:P\to X$, and whose morphisms are cartesian squares
\begin{equation}\label{torsors}
\begin{CD}
P'@>{g}>> P\\
@VVV @VVV\\
T'@>{h}>> T
\end{CD}
\end{equation}
of torsors such that the $G$-equivariant map $f':P'\to X$ is $g$ composed with the $G$-equivariant map $f:P\to X$. An object of $D^b_c([X/G])$ is the data of an $F_T\in D^b_{c}(T)$ for any torsor $P\to T$ as above, together with an isomorphism $h^{-1}(F_T)\to F_{T'}$ for each square (\ref{torsors}); these isomorphisms have to satisfy the obvious compatibility condition, cf. \cite[2.4.3]{berlunts}. A cohomology class $\in H^*([X/G],\Q_\ell)$ can be seen as the data of a class $\alpha_T\in H^*(T,\Q_\ell)$ for each couple $(P\to T,P\to X)$ as above, such that $h^*(\alpha_T)=\alpha_{T'}$ for each square (\ref{torsors}).

First of all, let us show that

\begin{lemma}\label{mainisoforstacks} $H^*(\Pi_{\dd,n}\setminus\Sigma_{\dd,n},\Q_\ell)$ is isomorphic to the tensor product of $H^*(\GL_{n+1},\Q_\ell)$ and the cohomology of the quotient stack $[(\Pi_{\dd,n}\setminus\Sigma_{\dd,n})/\GL_{n+1}(K)]$.
\end{lemma}
%{\bf Proof.} %In lemma \ref{explcoeffscharp} we constructed certain cohomology classes $\aaa^{\dd,n}_i(K,A)\in H^{2i-1}(\Pi_{\dd,n,K}\setminus\Sigma_{(d),n,K},A),i=1,\ldots,n+1,A=\Z/\ell^r$ such that their pullbacks under any orbit map are $\m_i^{\dd,n}\cc_i^{n+1}(K,A)$ where $\cc_i^{n+1}(K,A)$ are the canonical generators of $H^*(\GL_{n+1}(K),A)$. Taking the limit as $r\to \infty$ and tensoring with $\Q_\ell$ we get classes $\aaa^{\dd,n}_i(K,\Q_\ell)\in H^{2i-1}(\Pi_{\dd,n,K}\setminus\Sigma_{(d),n,K},\Q_\ell)$

{\bf Proof.} Let $p:P\to T$ be a $\GL_{n+1}(K)$-torsor and let $f:P\to \Pi_{\dd,n}\setminus\Sigma_{\dd,n}$ be a $\GL_{n+1}(K)$-equivariant map. The classes $\aaa^{\dd,n}_i(K,\Q_\ell)\in H^{2i-1}(\Pi_{\dd,n,K}\setminus\Sigma_{(d),n,K},\Q_\ell)$ from lemma \ref{explcoeffscharp} pulled back to $P$ using $f$ give a map in $D^b_c(T)$ from a direct sum of shifts of the constant sheaf to $Rp_* \underline{\Q_\ell}_P$.

This map is an isomorphism. Indeed, take a closed point $i:Spec (K)\to T$ of $T$;
let $p':P'\to Spec (K)$ be the pullback of $p$ under $i$ and let $\bar\imath$ be the corresponding map $P'\to P$, so that we get a Cartesian diagram
$$
\begin{CD}
P'@>{\bar\imath}>> P\\
@V{p'}VV @VV{p}V\\
Spec(K)@>{i}>> T
\end{CD}
$$
Since every torsor is locally trivial in the \'etale topology, the base change map of functors  $i^{-1}Rp_*\to Rp'_* \bar\imath^{-1}$ is an isomorphism, and so the above map $\sum \underline{\Q_\ell}_T[-a_i]\to Rp_* \underline{\Q_\ell}_P$ is an isomorphism as well. (Here $a_i$'s are the degrees of the elements of some homogeneous $\Q_\ell$-basis of $H^*(\GL_{n+1}(K),\Q_\ell)$.)

Now take $P$ to be $V_{n+1}(K^r)\times (\Pi_{\dd,n}\setminus\Sigma_{\dd,n})$ where $V_{n+1}(K^r)$ is the Stiefel manifold of $n+1$-frames in $K^r$. We have a natural $\GL_{n+1}(K)$-equivariant map $f:P\to\Pi_{\dd,n}\setminus\Sigma_{\dd,n}$ and a commutative diagram
$$
\begin{CD}
P@>{f}>> \Pi_{\dd,n}\setminus\Sigma_{\dd,n}\\
@V{p}VV @VV{q}V\\
[P/\GL_{n+1}(K)]@>>> [(\Pi_{\dd,n}\setminus\Sigma_{\dd,n})/\GL_{n+1}(K)]
\end{CD}
$$

For any segment $[a,b]\subset\Z$ of length $\leq r$ the bottom arrow induces an equivalence of categories $D^I_{\GL_{n+1}(K),c}(\Pi_{\dd,n}\setminus\Sigma_{\dd,n})\to D^I_{\GL_{n+1}(K),c}(P)$ where $D^I$ stands for $D^{\leq b}\cap D^{\geq a}$, cf. \cite[2.2.2, Corollary]{berlunts}. The quotient $P/\GL_{n+1}(K)$ is an algebraic variety (it is the Grassmannian of $n+1$-planes in $K^r$) and $D^b_{\GL_{n+1}(K),c}(P)$ is equivalent to $D^b (P/\GL_{n+1}(K))$, cf. \cite[Theorem 2.6.2]{berlunts}. Let $p:P\to P/\GL_{n+1}(K)$ be the quotient map.

As we saw above, $Rp_*\underline{\Q_\ell}_P=\sum \underline{\Q_\ell}_{P/\GL_{n+1}(K)}[-a_i]$. On the other hand, $Rf_*\underline{\Q_\ell}_P$ is $\underline{\mathbb{Q}_\ell}_{\Pi_{\dd,n}\setminus\Sigma_{\dd,n}}$ plus a sum of sheaves in degrees $\geq r$. So for $a=\left[\frac{r}{2}\right]$
the images in $D^{[-a,a]}_{\GL_{n+1}(K),c}(\Pi_{\dd,n}\setminus\Sigma_{\dd,n})$ of
%we have
%$$\tau^{\leq a}\tau^{\geq -a}Rq_* \underline{\Q_\ell}_{\Pi_{\dd,n}\setminus\Sigma_{\dd,n}}=\tau^{\leq a}\tau^{\geq -a}\underline{\Q_\ell}_{[(\Pi_{\dd,n}\setminus\Sigma_{\dd,n})/\GL_{n+1}(K)] }$$ 
$Rq_* \underline{\Q_\ell}_{\Pi_{\dd,n}\setminus\Sigma_{\dd,n}}$ and $\sum \underline{\Q_\ell}_{[(\Pi_{\dd,n}\setminus\Sigma_{\dd,n})/\GL_{n+1}(K)]}[-a_i]$ coincide. 
%as objects of $D^{[-a,a]}_{\GL_{n+1}(K),c}(\Pi_{\dd,n}\setminus\Sigma_{\dd,n})$.
Taking $r$ large enough, we see that $$Rq_* \underline{\Q_\ell}_{\Pi_{\dd,n}\setminus\Sigma_{\dd,n}}\cong \sum \underline{\Q_\ell}_{[(\Pi_{\dd,n}\setminus\Sigma_{\dd,n})/\GL_{n+1}(K)]}[-a_i].$$
The lemma is proven. $\clubsuit$

\begin{flushright}
{\sc Alexey Gorinov\\
Department of Mathemats\\
National Research University ``Higher School of Economics''\\
Moscow\\
Russia}\\
\url{agorinov@hse.ru, 
gorinov@mccme.ru}
\end{flushright}

\end{document}